\documentclass[12pt,reqno]{amsart}
\usepackage{pgfplots}
\usepackage{tikz}
\pgfplotsset{compat=newest}
\pgfplotsset{ytick style={draw=none}}
\pgfplotsset{xtick style={draw=none}}
\usepackage{amsmath,amsfonts,amssymb,amsthm}
\usepackage{graphicx}
\graphicspath{ }
\usepackage{mathrsfs}
\usepackage{epstopdf}
\usepackage{csquotes}
\usepackage{wrapfig}
\usepackage{accents}
\usepackage{caption}
\usepackage{subcaption}
\usepackage{calligra}
\usepackage[colorlinks]{hyperref}
\usepackage{kantlipsum}
\hypersetup{colorlinks, 
	breaklinks,
	linkcolor=red,
	citecolor=red,
	linktocpage=true}

\numberwithin{figure}{section}

\theoremstyle{plain}
\newtheorem{thm}{Theorem}[section]
\newtheorem{lem}[thm]{Lemma}
\newtheorem{prop}[thm]{Proposition}

\theoremstyle{definition}
\newtheorem{defn}{Definition}[section]

\newtheorem{exmp}{Example}[section]
\theoremstyle{remark}
\newtheorem*{rem}{Remark}

\usepackage{mathtools}
\makeatletter
\@namedef{subjclassname@2020}{%
	\textup{2020} Mathematics Subject Classification}

\makeatother
\author[M. A. Sarkar]{MABUD ALI SARKAR$^{*}$}
\address{Department of Mathematics\\ The University of Burdwan \\ Burdwan-713101, India.}
\email{mabudji@gmail.com$^{*}$}
\author[A. A. Shaikh]{Absos Ali Shaikh$^{\dagger}$}
\address{Department of Mathematics\\ The University of Burdwan \\ Burdwan-713101, India.}
\email{aashaikh@math.buruniv.ac.in$^{\dagger}$}

\title{On the bases of the image of $2$-adic logarithm on the group of principal units}

\begin{document}
	
		\begin{abstract}  This paper computes the bases of the image of $2$-adic logarithm on the group of the principal units in all 7 quadratic extensions of $\mathbb{Q}_2$. This helps one to understand the free module structure of $2$-adic logarithm at arbitrary points on its domain. We discuss some applications at the end.
		\end{abstract} 

		\maketitle

		\section*{Introduction}
		The $p$-adic logarithm $\log_p(1+x)=\sum_{n\geq 1}(-1)^{n-1}\frac{x^n}{n} \in K[[x]]$ converges for all $x \in \mathfrak{m}_K$. Following Iwasawa \cite{WA}, the $p$-adic logarithm extends to $p^{\mathbb{Z}}(1+\mathfrak{m}_{\mathbb{C}_p}) \subset \mathbb{C}_p^{\times}$ such that $\log_p(p)=0$. It is known that the $p$-adic logarithm defines an isomorphism from $1+\mathfrak{m}_K^r \to \mathfrak{m}_K^r$ if $r> \frac{e}{p-1}$, where $e$ is the ramification index. However, the case $r=1$ is different i.e, $\log_p: 1+\mathfrak{m}_K \to \mathfrak{m}_K$ is not an isomorphism. The image $\log_p(1+\mathfrak{m}_K)$ is not known for arbitrary extension $K$ of $\mathbb{Q}_p$. 
		
		In $p$-adic number theory, the $p$-adic logarithm plays an essential role. From (\cite{BT}, \cite{KI1}, \cite{KI}, \cite{JC}, \cite{WA}), we know the image of $p$-adic logarithm on the groups of principal units is very crucial in the theory of Iwasawa. More specifically in \cite{KI}, Iwasawa explicitly computed the formulas for the norm residue symbol by using the image of $p$-adic logarithm on the group of principal units of certain cyclotomic extension of $\mathbb{Q}_p$, which is our main motivation. There are plenty of other potential applications of the image/value of $p$-adic logarithm. In \cite{GG}, the formula for normalized $p$-adic regulator contains the image of $p$-adic logarithm on the principal units. In the study of $p$-adic $L$-functions there appears $p$-adic logarithm. To compute the $p$-adic $L$-function $L(1, \chi)$, one of the key challenges is the computation of $p$-adic logarithm of arbitrary elements in its domain, e.g., \cite{FY}.  The value of $p$-adic logarithm has an important role in the formula of Sen operator (\cite{LB}) in the theory of $p$-adic representation as well as in the formula of $p$-adic heights \cite{BM} in the theory of arithmetic dynamics. More applications about the image of $p$-adic logarithm are discussed in Section \ref{s3}. Therefore, it is of novelty to understand the image/value of $p$-adic logarithm on its domain.  Since by Proposition \ref{p1.4}, $\log_p(1+\mathfrak{m}_K)$ is a $\mathbb{Z}_p$-module, once we know the bases of $\log_p(1+\mathfrak{m}_K)$, we can understand the value of the $p$-adic logarithm at an arbitrary point in its domain by forming a linear combination of that basis elements with scalars from $\mathbb{Z}_p$. Thus, this current paper undertakes an initial challenge towards computing the image of $p$-adic logarithm on the group of principal units, by computing the bases of $p$-adic logarithm. A basis of a free module is as important as a bases of a vector space. In Theorem \ref{t2.13} and Theorem \ref{t2.14} of Section \ref{s2}, we compute the basis of $2$-adic logarithm for all $7$ quadratic extensions of $\mathbb{Q}_2$. Moreover, we find the exact images of $2$-adic logarithm in some cases. In Section \ref{s3}, we discuss some more applications. \\
		\subsection*{Notations} $\mathbb{Q}_p$ is the p-adic number field with ring of integers $\mathbb{Z}_p$, unique maximal ideal $p \mathbb{Z}_p$ and residue field $\mathbb{F}_p:=\mathbb{Z}_p/p \mathbb{Z}_p$. Let $\bar{\mathbb{Q}}_p$ be the algebraic closure and take the $p$-adic completion $\mathbb{C}_p:=\widehat{\bar{Q}}_p$ with maximal ideal $\mathfrak{m}_{\mathbb{C}_p}$ and the units $\mathbb{C}_p^{\times}$. Let $K$ be a finite extension of $\mathbb{Q}_p$ with ring of integers $\mathcal{O}_K$, unique maximal ideal $\mathfrak{m}_K$. The quotient field $\kappa_K:=\mathcal{O}_K/\mathfrak{m}_K$ is the residue field of $\mathcal{O}_K$. 
		We also have the standard $p$-adic additive valuation $v_p: \mathbb{Q}_p \to \mathbb{Z} \cup \{ \infty\}$.
		We extend the valuation $v_p$ on $\mathbb{Q}_p$ to the finite extension $K$ denoted by $v$ that makes $v(z)=v_p(z)$ whenever $z \in \mathbb{Q}_p$ by the following formula $v(\lambda)=\frac{1}{[K:\mathbb{Q}_p]} v_p \left(\mathcal{N}_{\mathbb{Q}_p}^K(\lambda)\right),$
		where $\mathcal{N}_{\mathbb{Q}_p}^K$ is the field-theoretic norm, the multiplicative mapping from $K^{\times}$ to $\mathbb{Q}_p$. We use the common notation $v$ throughout the paper.

		\section{Formal group law and $p$-adic logarithm}
		In this section, we construct the $p$-adic logarithm from $1$-dimensional multiplicative formal group law.
		\begin{defn} \cite{M}
			Let $R$ be a commutative ring. A formal group law over $R$ is a formal power series $F(X,Y) \in R[[X,Y]]$ with the following properties
			\begin{enumerate}
				\item[(i)]  $F(X,Y)=X+Y+ \text{higher degree terms}$,
				\item[(ii)]  $F(X,F(Y,Z))=F(F(X,Y),Z),$
				\item[(iii)] $F(X,0)=X$ and $F(0,Y)=Y$,
				\item[(iv)]  There exists a unique power series $\iota (T) \in R[[T]] $ such that $$F(X, \iota(X))=F(i(Y), Y)=0,$$which is the inverse power series in the formal group law. 
			\end{enumerate}
			If in addition $F(X, Y)=F(Y, X)$ holds then $F$ is called commutative. The commutativity property comes trivially as long as $R$ has no elements that are both torsion and nilpotent.
		\end{defn}
		\begin{exmp}
			The $1$-dimensional formal additive group law is defined by $F(X, Y)=X+Y$ and the $1$-dimensional formal multiplicative group law is defined by $\mathcal{M}(X, Y)=X+Y+XY.$
		\end{exmp} 
		Now, we discuss maps between formal group laws.
		\begin{defn} \label{def:2.2} \cite{J}
			Let $F,~G$ be two formal group laws over a ring $R$.  A homomorphism $f: F \to G$ over $R$ is a power series $f(T) \in R[[T]]$, having no constant term, satisfying $$f(F(X,Y))=G(f(X), f(Y)).$$ If further there exists another homomorphism $\psi \in R[[T]]$, $\psi: G \to F$ such that
			$$f(\psi(T))=\psi(f(T))=T,$$ then we say that $F$ and $G$ are isomorphic over $R$ while $f$ is an isomorphism over $R$. An endomorphism $f: F \to F$ satisfies $f(F(X,Y))=F(f(X),f(Y)).$
		\end{defn} 
		\begin{defn} \cite{J}
			An invariant differential of a formal group law $F$ over an arbitrary ring $R$ is a differential form $ \omega(T)=p(T)dT \in R [[T]]dT$ that satisfies $p(F(T,~S))F_X(T,~S)=p(T)$, where $F_X(X,~Y)$ is the partial derivative with respect to its first variable. An invariant differential is said to be normalized if $p(0)=1$.
			Then there exists a unique normalized invariant differential of $F$ given by the formula $$\omega(T)=F_X(0,T)^{-1}dT.$$
		\end{defn}
		We now construct the formal logarithm, of a formal group law, obtained by integrating the invariant differential. But integrating an invariant differential form $\omega(T)$ produces denominators in each coefficient. So we assume $R$ to be a commutative ring of characteristic $0$ so that we don't have further difficulty of division. The following construction of the $p$-adic logarithm can be found in \cite[Section~5.8]{M}, but for the convenience of readers we add a proof here, providing minor details:   
		\begin{prop} Assume that $R$ be a torsion-free ring and let $R'=R \otimes_{\mathbb{Z}} \mathbb{Q}$. If $\omega(T)=p(T)dT$ is a normalized invariant differential of a formal group law $F$. Then the formal group law $F$ induces a formal logarithmic series $\log_F$ on $R'$ given by $$\log_F(T)=\int_0^T \omega(\tau)=T+\frac{a_2}{2}T^2+\frac{a_3}{3}T^3+\cdots \in R'[[T]],$$ satisfying $\log_F(F(X,Y))=\log_F(X)+\log_F(Y)$.
		\end{prop}   
		\begin{proof}
			Define  the map $R \to R'=R \otimes_{\mathbb{Z}} \mathbb{Q}$ by $ r \mapsto r \otimes 1$ for $r \in R$. The kernel of this map is the set of torsion elements of $R$. But as $R$ is the torsion-free ring, the kernel is trivial. Therefore there is a one-to-one crorrespondence between $R$ and $R \otimes_{\mathbb{Z}} \mathbb{Q}$. This allows us to take coefficients of the formal logarithm in $\mathbb{Q}$. 
			
			Now as $p(T)dT$ is normalised invariant differential, we have $p(0)=1$ and so we can assume $$p(T)=F_X(0,T)^{-1}=1+ \sum_{i=1}^{\infty}a_iT^i \in R'[[T]].$$
			Integrating, we get $$ \log_F(T)=\int_0^T \omega(\tau)=\int_0^T \frac{1}{F_X(0,T)}dT=T+\frac{a_1}{2}T^2+\frac{a_2}{3}T^3+\cdots \in R'[[T]].$$
			Now it is remains to show that $\log_F(F(X,Y))=\log_F(X)+\log_F(Y).$ To do this we take partial derivatives with respect to $X$ in the associativity relation $F(X,F(Y,Z))=F(F(X,Y),Z)$, we get
			\begin{equation} \label{1}
			F_X(X,F(Y,Z))=F_X(F(X,Y),Z) \cdot F_X(X,Y).
			\end{equation}
			Let $p(X)$ be the formal power series with coefficients in $R$ defined by
			\begin{equation} \label{2}
			p(X) \cdot F_X(0,X)=p(0)=1.
			\end{equation}
			Now $\log_F(X)$ is the unique formal power series with coefficients in $R'=R \otimes_{\mathbb{Z}} \mathbb{Q}$  such that
			\begin{equation} \label{3}
			\frac{\partial}{\partial X} \log_F(X)=p(X).
			\end{equation}
			Then using the changes $X=0,\ Y \to X, \ Z \to Y$ in \eqref{1} and using \eqref{2}, we get
			\begin{equation} \label{4}
			\frac{\partial}{\partial X}\log_F(F(X,Y))=\frac{\partial}{\partial X} \log_F(X).
			\end{equation}
			Integrating \eqref{4}, we get
			\begin{equation} \label{5}
			\log_F(F(X,Y))=\log_F(X)+C(Y), \ \text{where $C(Y)$ is an arbitrary functions in $Y$},
			\end{equation}
			which can be written as 
			\begin{equation} \label{6}
			\log_F(F(X,Y))-\log_F(X)=C(Y).
			\end{equation}
			Putting $X=0$ in \eqref{6}, we have 
			\begin{align*}
				&\log_F(F(0,Y))-\log_F(0)=C(Y) \\
				\Rightarrow &\log_F(Y)=C(Y), ~\text{since}~ \log_F(0)=0~\text{and}~ F(0,Y)=Y.
			\end{align*} 
			Hence from \eqref{6}, we conclude $\log_F(F(X,Y))=\log_F(X)+\log_F(Y).$
			This completes the proof.

		\end{proof}

		\begin{prop} \label{p1.2} Let $K$ be a finite extension of $\mathbb{Q}_p$ with ring of integers $\mathcal{O}_K$ and unique maximal ideal $\mathfrak{m}_K$. Then the
			1-dimensional multiplicaive formal group law $\mathcal{M}(x,y)=x+y+xy$ defined over $\mathcal{O}_K$ induces the $p$-adic logarithm (analogue to classical logarithm) given by $$\log_{\mathcal{M}}(t)=\log_p(1+t)=t-\frac{t^2}{2}+\frac{t^3}{3}-\cdots$$, and converges in $\mathfrak{m}_K$.
		\end{prop}
		\begin{proof}
			If $\pi \in \mathcal{O}_K$ is a uniformizer, we can choose either $\pi=\sqrt[p-1]{-p}$ or $\pi=\zeta_p-1$, for some $p^{th}$ root of unity.
			It can be checked that the differential form $\omega(x)=\frac{dx}{1+x}$ is the invariant differential for $\mathcal{M}(x,y)=x+y+xy$, and hence the corresponding logarithmic series $\log_{\mathcal{M}}$ is obtained by $$\log_{\mathcal{M}}(x)=\int_0^x \omega(t)dt=\int_0^x \frac{dt}{1+t}=x-\frac{x^2}{2}+\frac{x^3}{3}-\cdots,$$ which satisfies
			$ \log_{\mathcal{M}}(\mathcal{M}(x,y))=\log_{\mathcal{M}}(x)+\log_{\mathcal{M}}(y).$
			This logarithmic series is called the p-adic logarithm and is denoted by 
			\begin{equation} \label{e7}
			\log_{\mathcal{M}}(x)=\log_p(1+x)=x-\frac{x^2}{2}+\frac{x^3}{3}-\cdots \in \mathcal{O}_K[[x]].
			\end{equation} This series converges in $\mathfrak{m}_K$. For, we note that the Cauchy-Hadamard convergence formula holds in the $p$-adic field as well. 
			If  $f(x)=\sum_{n \geq 0} a_n x^n \in K[[x]]$, by Cauchy-Hadamard formula, the radius of convergence $r \in [0,\infty]$ is given by 
			$$ \frac{1}{r}=\varlimsup_{n \to \infty} \sqrt[n]{|a_n|_p}. $$
			Therefore, for $p$-adic logarithm \eqref{e7}, we have $a_n=(-1)^{n-1} \frac{1}{n}$. Then, we have
			$$ \frac{1}{r}=\varlimsup_{n \to \infty} \sqrt[n]{\left| \frac{(-1)^{n-1}}{n} \right|_p}=\varlimsup_{n \to \infty} \frac{1}{\sqrt[n]{|n|_p}}=\varlimsup_{n \to \infty}p^{v_p(n)/n}. $$
			Since $1 \leq p^{v_p(n)} \leq n$, we get $ 1 \leq p^{v_p(n)/n} \leq \sqrt[n]{n}$. Also $\sqrt[n]{n} \to 1$ as $n \to \infty$. In particular $ \sqrt[n]{|n|_p} \to 1$ as $n \to \infty$. Hence the radius of convergence $r=1$ i.e., $-1 <x \leq 1$. At the end point $x=1$, we have $a_n=\pm\frac{1}{n}$ and its $p$-adic norm $\frac{1}{|n|_p}$ doesn't tends to $0$ because it is $1$ when $p  \nmid n $. Therefore the disc of convergence of the $p$-adic logarithm in $K$ is the open unit disc $\{x \in K: |x|_p<1 \}=\mathfrak{m}_K$. 
		\end{proof}
		
		Thus the important thing is to note that the $1$-dimensional multiplicative formal group law $\mathcal{M}$ gives a group structure on $\mathfrak{m}_K$ as follows:
		\begin{lem} 
			Let $F$ be a formal group law over $\mathcal{O}_K$. Then the set $\mathfrak{m}_K$ becomes a group under the group law $F$, via the law of combination $\alpha+_F \beta=F(\alpha,\beta)$. \label{2.3} \end{lem}
		\begin{proof}
			The proof is obvious using the previous definitions and basic properties of formal group law. The only thing that needs to ensure is the convergence of the infinite series $F(\alpha, \beta)$, which follows from the completeness of $\mathfrak{m}_K$ under $\mathfrak{m}_K$-adic topology. 
		\end{proof}
		The following standard proposition describes the module structure of $p$-adic logarithm on principal units. The proof is based on basic concept of $p$-adic logarithm, e.g., \cite[Chap.~IV]{NK}, \cite[Sec.~4,~Chap.~5]{AR}. 
		\begin{prop} \label{p1.4}
			$\log_p(1+\mathfrak{m}_K)=\log_{\mathcal{M}}(\mathfrak{m}_K)$ is a $\mathbb{Z}_p$-module.
		\end{prop}
		\begin{proof}
			The action of the ring $\mathbb{Z}_p$ on $\log_{p}(1+\mathfrak{m}_K)$ can be described as $\mathbb{Z}_p \times \log_{p}(1+\mathfrak{m}_K) \to \log_{p}(1+\mathfrak{m}_K)$ defined by $a \cdot \log_p(b)=\log_p(b^a)$. This is valid because for $a \in \mathbb{Z}_p$, the binomial expension $$(1+b)^a=\sum_{n=0}^{\infty}\frac{a(a-1) \cdots (a-n+1)}{n!} b^n \ \text{converges}.$$ 
		\end{proof}
		\begin{rem}
			The notation $\log_{p}(1+x)$ of $p$-adic logarithm indicates the field $\mathbb{Q}_p$, we are considering. While the notation $\log_{\mathcal{M}}(x)$ indicates that the $p$-adic logarithm is induced by the $1$-dimensional multiplicative formal group law $\mathcal{M}(X, Y)=X+Y+XY$. The use of the formal group $\mathcal{M}$ makes notation and concept much clearer here, because, the original multiplicative group $1+\mathfrak{m}_K$ gets described, using $\mathcal{M}$, as just $\mathfrak{m}_K$. Most of the time, in this article we use the notation $\log_{\mathcal{M}}$ if no confusion arises. 
		\end{rem}
		\section{Compuation of basis of $2$-adic logarithm} \label{s2}
		In this section, we compute the bases of the images of the $2$-adic logarithm $\log_2(1+\mathfrak{m}_K)$, where $\mathfrak{m}_K$ are maximal ideals in the $7$ quadratic extensions $K=\mathbb{Q}_2(\sqrt{d})$ of $\mathbb{Q}_2$, $d=-1, \pm 2, \pm 3, \pm 6$. Throughout this section, we use the notation $\log_{\mathcal{M}}(x)$ instead of the notation $\log_2(1+x)$ to use the properties of the $1$-dimensional multiplicative formal group law $\mathcal{M}$.

		The multiplication-by-$n$ map $n \to [n]_F$ is an endomorphism of a formal group law $F$ from the ordinary integers $\mathbb{Z}$ into the ring of endomorphisms $\text{End}_{\mathcal{O}_K}(F)$ satisfying:
		\begin{enumerate}
			\item[(i)] $[0]_F=0$ and for $n \geq 0$, $[n+1]_F(X)=F([n]_F(X),X),$
			\item[(ii)] $[n]_F=\iota_F \circ [-n]_F$, if $n<0$.
		\end{enumerate}
		We have the following relations for the $1$-dimensional multiplicative formal group law $\mathcal{M}$:
		\begin{lem} \label{l2.1}
			We have the following consequences,
			\begin{enumerate}
				\item [(i)] $\mathcal{M}(-2,-2)=[2]_{\mathcal{M}}(-2)=0$.
				\item [(ii)] $\log_{\mathcal{M}}(-2)=0$.
				\item [(iii)]$\log_{\mathcal{M}}(2)=0 \ (\text{mod} \ 4 \mathcal{O}_K)$.
				\item [(iv)]  $\mathcal{M}(2,2)=[2]_{\mathcal{M}}(2)=0 \ (\text{mod} \ 4 \mathcal{O}_K)$.
			\end{enumerate}
		\end{lem}
		\begin{proof}
			\begin{enumerate}
				\item[(i)] We have $\mathcal{M}(x,y)=x+y+xy$. Then, $\mathcal{M}(-2,-2)=-2-2+4=0.$
				\item[(ii)] It follows from $(i)$.
				\item[(iii)] $\log_{\mathcal{M}}(2)=\mathcal{M}(2,2)=2+2+4=8 \equiv 0 \ (\text{mod} \ 4 \mathcal{O}_K)$.
				\item[(iv)] It follows from $(iii)$.
			\end{enumerate}
		\end{proof}
		
		We note the following result:
		\begin{thm} \cite{GS}
			There are $7$ quadratic extensions of $\mathbb{Q}_2$ given by        $\mathbb{Q}_2(\sqrt{d})$ for $d=-1,\pm 2,\pm 3, \pm 6$.  All these extensions except $\mathbb{Q}_2(\sqrt{-3})$ are ramified.
		\end{thm}
		Each of these $7$ extensions has multiple equally-good descriptions, for instance $\mathbb{Q}_2(\sqrt{-3})=\mathbb{Q}_2(\sqrt{5})$. Of these seven, $\mathbb{Q}_2(\sqrt{-3})$ is the only unramified extension and $\mathbb{Q}_2(\sqrt{-1})$ is the extension with additional $2$-power roots of unity.
		\begin{lem} \label{l4.2}
			Let $K$ be a finite extension of $\mathbb{Q}_2$. Then the $\log_{\mathcal{M}}$ series is convergent on the maximal ideal $\mathfrak{m}_K$ and maps the set $2 \mathfrak{m}_K$ onto itself in a one-to-one way.
			$ \label{eq: 4.1}$ \end{lem}
		\begin{proof}
			The Proposition \ref{p1.2} shows the $p$-adic logarithm is convergent in the maximal ideal.
			For the second claim, looking at the power series $$\log_{\mathcal{M}}(x)=\frac{1}{2} \log_{\mathcal{M}}(2x)=x-x^2+\frac{4x^3}{3}-2x^4+\frac{16x^5}{5}-\frac{16x^6}{3}+ \cdots ,$$
			whose coefficients are all in $\mathbb{Z}_2$. Since $\log_{\mathcal{M}}$ has its first-degree coefficient in $\mathbb{Z}_2^{\times}$, the units of $\mathbb{Z}_2$, it has an inverse power series, $\log_{\mathcal{M}}^{-1}=x+x^2+\cdots \in \mathbb{Z}_2[[x]]$, of course we know that this inverse function is $\frac{1}{2} (exp(2x)-1)$. This inverse function is also convergent for $x \in \mathfrak{m}_K$. Thus we have two series in $\mathbb{Z}_2[[x]]$, inverse to each other and thus $\log_{\mathcal{M}}$ is one-to-one and onto $2 \mathfrak{m}_K$.
		\end{proof}
		\begin{lem} \label{l2.4}
			The image of $4 \mathbb{Z}_2$ under 2-adic logarithm $\log_{\mathcal{M}}$ is $4 \mathbb{Z}_2$. Furthermore, if $z=2^m u$ for a 2-adic unit $u$ and with $m \geq 2$, then $\log_{\mathcal{M}}(z)=2^mu'$ for a 2-adic unit $u'$. In other words, if $v_2(z) \geq 2$, then $v_2(\log_{\mathcal{M}}(z))=v_2(z)$.
		\end{lem} 
		\begin{proof}
			Applying Lemma $\ref{eq: 4.1}$  the proof follows. Therefore $\log_{\mathcal{M}}$ maps $2 \cdot 2 \mathbb{Z}_2$ onto itself in a one-to-one way.
		\end{proof}
		\begin{lem} \label{eq:4.3}
			For $p>2$, the $\mathbb{Q}_p$-series $\frac{1}{p} \log_{\mathcal{M}}(px)$ has all coefficients in $\mathbb{Z}_p$ and all but first degree coefficients are in $p \mathbb{Z}_p$.
		\end{lem}
		\begin{proof}
			Writing out the series $\frac{1}{p}\log_{\mathcal{M}}(px)=x-\frac{px^2}{2}+\frac{p^2x^3}{3}-\cdots,$
			and looking at the denominator, the $n$-th denominator in the coefficient satisfy $v_p(n) \leq \log_p(n)$, while the numerator has $v_p(p^{n-1})=n-1$. Since for $p>2$ and $n>1$, we have the inequality $n>\log_p(n)+1$, the result follows immediately. 
		\end{proof}
		\begin{lem}
			Let $p>2$, the image of $p \mathbb{Z}_p$ under the $p$-adic logarithm is $p \mathbb{Z}_p$. Furthermore, if $z=p^mu$ for a $p$-adic unit $u$ with $m \geq 1$, then $\log_{\mathcal{M}}(z)=p^mu'$ for a $p$-adic unit $u'$. In other words, if $v_p(z) \geq 1$, then $v_p(\log_{\mathcal{M}}(z))=v_p(z)$.
		\end{lem}
		\begin{proof}
			Applying previous Lemma $ \ref{eq:4.3}$, the result follows immediately.
		\end{proof}
		Let $C_n$ refer to the cyclic group of order $n$. We can think of it as the additive group $\mathbb{Z}/n \mathbb{Z}$, but in many cases, the group will be written multiplicatively.  
		\begin{lem} \label{l2.7}
			Let $K$ be a quadratic extension of $\mathbb{Q}_2$, its ring of integers $\mathcal{O}_K$, and the maximal ideal $\mathfrak{m}_K=\pi \mathcal{O}_K$ for uniformizer $\pi \in \mathcal{O}_K$. Then, in the additive group $\mathcal{O}_K$, the subgroup $\mathcal{O}_K/2 \mathcal{O}_K$ has order $4$, and the factor group $\mathcal{O}_K/2 \mathcal{O}_K$ has the structure $C_2 \oplus C_2$. The same is true for the quotient $\mathfrak{m}_K/2 \mathfrak{m}_K$. Under the formal group addition $+_{\mathcal{M}}$, the subgroup $2 \mathfrak{m}_K$ of $\mathfrak{m}_K$ still has index $4$.
		\end{lem}
		\begin{proof}
			Using the $\text{structure theorem of a finitely generated module over a principal ideal domain}$, we see that a finitely generated module over $\mathbb{Z}_2$ will have the decomposition $\mathbb{Z}_2^m \oplus \mathcal{F}$, where the first part $\mathbb{Z}_2^m$ is the free part, which is torsion-free and  $\mathcal{F}$ is a direct sum of modules of the form $\mathbb{Z}_2/2^m \mathbb{Z}_2$, which is torsion group. Thus $\mathcal{O}_K \cong \mathbb{Z}_2^{m}$ for some $m$, the maximal cardinality of a $\mathbb{Z}_2$- linearly independent subset of $\mathcal{O}_K$. But such a set is also a basis of the quadratic extension $K$ over $\mathbb{Q}_2$. \\ Thus we have $\mathcal{O}_K/2 \mathcal{O}_K \cong \mathbb{Z}_2^2/2 \mathbb{Z}_2^2 \cong (\mathbb{Z}_2/2 \mathbb{Z}_2)^2$, the second isomorphism follows from a map $f: \mathbb{Z}_2^2=\mathbb{Z}_2 \times \mathbb{Z}_2 \to (\mathbb{Z}_2/2\mathbb{Z}_2)^2=\mathbb{Z}_2/2\mathbb{Z}_2 \times \mathbb{Z}_2/2 \mathbb{Z}_2$ defined by $f(a,b)=(a+2 \mathbb{Z}_2, b+2 \mathbb{Z}_2)$ for $(a,b) \in \mathbb{Z}_2 \times \mathbb{Z}_2$. This map is a surjective homomorphism. The kernel is $2 \mathbb{Z}_2^2$ because $f(a,b)=(0+2\mathbb{Z}_2,0+2\mathbb{Z}_2) \Leftrightarrow (a,b)=(2u,2v), \ u,v \in \mathbb{Z}_2 \Leftrightarrow (a,b)=2(u,v) \in 2(\mathbb{Z}_2 \times \mathbb{Z}_2)=2\mathbb{Z}_2^2$.  Also $C_2^2=(\mathbb{Z}_2/2 \mathbb{Z}_2)^2$ and hence $\mathcal{O}_K/2 \mathcal{O}_K \cong C_2^2 \cong C_2 \oplus C_2$. This proves the first assertion. Since $\mathfrak{m}_K=\pi \mathcal{O}_K$, the multiplication-by-$\pi$ map is a $\mathbb{Z}_2$-homomorphism from $\mathcal{O}_K$ onto $\pi \mathcal{O}_K$, and the next assertion follows i.e., $\mathfrak{m}_K/2 \mathfrak{m}_K \cong C_2 \oplus C_2$. \\ The structure of $\mathfrak{m}_K$ under the addition $+_{\mathcal{M}}$ may be different from the familiar additive structure, but at least $2 \mathfrak{m}_K$ is still a subgroup of $\mathfrak{m}_K$ of index four. For, if $K$ is ramified extension over $\mathbb{Q}_2$, then $\mathfrak{m}_K=\pi \mathcal{O}_K$, $\mathfrak{m}_K^2=\pi^2 \mathcal{O}_K$, and $\mathfrak{m}_K^3=\pi \cdot \pi^2 \mathcal{O}_K=\pi \cdot 2 \mathcal{O}_K$, since $\pi^2 \mathcal{O}_K=2 \mathcal{O}_K$. Then $[\pi^i \mathcal{O}_K : \pi^{i+1} \mathcal{O}_K]=2,$ for $i \geq 0$, because $\mathcal{O}_K/\pi \mathcal{O}_K=\mathcal{O}_K/\mathfrak{m}_K=\kappa_K=\mathbb{F}_2$.  \\ If $K$ is an unramified extension of $\mathbb{Q}_2$, we take the basis $\{1, \omega \}$ of $\mathcal{O}_K$ over $\mathbb{Z}_2$, where $\omega$ is the primitive cube root of unity, i.e., $1+\omega+\omega^2=0$. Then, $$2 +_\mathcal{M} 2 \omega=2+2 \omega+4 \omega \equiv 2+2 \omega \ ( \text{mod}  \ 4 \mathcal{O}_K).$$ 
			We note, 
			
			\[
			\left\{\begin{aligned}
			For, & \ z=2, \ 2+_{\mathcal{M}}2=\mathcal{M}(2,2)\equiv 0 \ (\text{mod} \ 4 \mathcal{O}_K), \\ For, & \ z=2 \omega, \ 2 \omega+_{\mathcal{M}} 2 \omega=\mathcal{M}(2 \omega, \omega) \equiv 0 \ (\text{mod} \ 4 \mathcal{O}_K), \\ For, & \ z=2+2 \omega, \ (2+2 \omega)+_{\mathcal{M}} (2+2 \omega)=\mathcal{M}(2+2\omega,2+2 \omega) \equiv 0 \ (\text{mod} \ 4 \mathcal{O}_K).
			\end{aligned} \right. 
			\]
			Thus we get $z+_{\mathcal{M}} z \equiv 0 \ (\text{mod} \ 4 \mathcal{O}_K)$ and hence the factor group $\mathcal{O}_K/2 \mathcal{O}_K$ has the structure $C_2 \oplus C_2$. 
		\end{proof} 
		\subsection{Main results}
		\begin{lem} \label{l2.8}
			If we consider $K=\mathbb{Q}_2,~\mathcal{O}_K=\mathbb{Z}_2,~\mathfrak{m}_K=2 \mathbb{Z}_2,~\pi=2$, then $\log_{\mathcal{M}}( \mathfrak{m}_K)=\log_{\mathcal{M}}(2 \mathfrak{m}_K)=2 \mathfrak{m}_K$.
		\end{lem}
		\begin{proof}
			We must show that if $z \in 2 \mathbb{Z}_2$, then $\log_{\mathcal{M}}(z) \in 4 \mathbb{Z}_2$; and that if $w \in 4 \mathbb{Z}_2$, then there is $z \in 2 \mathbb{Z}_2$ with $\log_{\mathcal{M}}(z)=w$.
			
			We start with $z \in 2 \mathbb{Z}_2$. We know that if $z \in 4 \mathbb{Z}_2$, then $\log_{\mathcal{M}}(z) \in 4 \mathbb{Z}_2$, by Lemma \ref{l2.4}, so that we need to consider only the case where $z \equiv 0~(\text{mod}~2)$ but $z \not\equiv 0~(\text{mod}~4)$: in other words, $z \equiv 2~(\text{mod}~4)$. Now consider $w=z+_{\mathcal{M}}(-2)=z-2-2z \equiv z-2 \equiv 0~(\text{mod} \ 4)$. We have $\log_{\mathcal{M}}(w) \in 4 \mathbb{Z}_2$ by Lemma \ref{l2.1}, and at the same time $\log_{\mathcal{M}}(w)=\log_{\mathcal{M}}(z)+\log_{\mathcal{M}}(-2)=\log_{\mathcal{M}}(z)$ because $\log_{\mathcal{M}}(-2)=0$. So $\log_{\mathcal{M}}(z) \in 4 \mathbb{Z}_2$. 
			
			The other direction is trivial: if $w \in 4 \mathbb{Z}_2$, then $w=\log_{\mathcal{M}}(z)$ for $z \in \mathbb{Z}_2$, and since $4 \mathbb{Z}_2 \subset 2 \mathbb{Z}_2$, we are done.
		\end{proof}
		\begin{proof}[\protect{Alternative proof of Lemma~\ref{l2.8}}]
			What is to be proved is that $\log_{\mathcal{M}}(2 \mathbb{Z}_2)=\log_{\mathcal{M}} (4 \mathbb{Z}_2)$. What the proof depends on is the fact that $\log_{\mathcal{M}}(2)=4u$ for $u$ a $\mathbb{Z}_2$-unit, in other words, that $v_2(\log_{\mathcal{M}}(2))=2$. For this, we do a direct computation using the series expansion for $\log_{\mathcal{M}}$: the first two terms, $2$ and $-2^2/2$, cancel each other, and the most interesting one is $-2^4/4=-4$. All other terms have $v_2$-value at least $3$, so the claim about $v_2(\log_{\mathcal{M}}(2))$ is verified.
			
			Now suppose $v(z)=1$, then $z+_{\mathcal{M}}2=z+2+2z$, in which $z \equiv 2 ~(\text{mod} \ 4)$ so that $v(z+2) \equiv 0~(\text{mod} \ 4)$, and hence it follows that $v(z+_{\mathcal{M}} 2) \geq 2$. Now, we have $\log_{\mathcal{M}}(z)=\log_{\mathcal{M}}(z+_{\mathcal{M}}2)-\log_{\mathcal{M}}(2)$, in which both elements on the right-hand side have $v$-value $2$, so that their difference on the left has $v$-value $\geq 2$.
		\end{proof}
		
		Now, in case $K$ is quadratic ramified over $\mathbb{Q}_2$, with ring of integers $\mathcal{O}_K$ and maximal ideal $\mathfrak{m}_K=\pi \mathcal{O}_K$, we have the chain of multiplicative subgroups of $1+\mathfrak{m}_K$:
		$$\cdots \subset 1+2 \mathfrak{m}_K=1+\pi^2 \mathfrak{m}_K=1+\pi^3 \mathcal{O}_K \subset 1+\pi \mathfrak{m}_K=1+\pi^2 \mathcal{O}_K \subset 1+\mathfrak{m}_K=1+\pi \mathcal{O}_K,$$
		where the inclusions are strict, and each index $[\pi^m \mathcal{O}_K:\pi^{m+1} \mathcal{O}_K ]=2$. Translating to the group structure of $\mathfrak{m}_K$ furnished by $+_{\mathcal{M}}$, we see the same inclusions: 
		$$2 \mathfrak{m}_K=\pi^2 \mathfrak{m}_K=\pi^3 \mathcal{O}_K \subset \pi \mathfrak{m}_K=\pi^2 \mathcal{O}_K \subset \mathfrak{m}_K=\pi \mathcal{O}_K.$$
		Now we will associate here $\log_{\mathcal{M}}(\pi)$ with $\log_{\mathcal{M}}(\mathfrak{m}_K)$. We know that $\log_{\mathcal{M}} (2 \mathfrak{m}_K)=2 \mathfrak{m}_K$, and since in the $+_{\mathcal{M}}$-group structure a representative of the nonzero element of $\pi \mathfrak{m}_K/2 \mathfrak{m}_K$ is $-2$, which has the property that $\log_{\mathcal{M}} (-2)=0$, we can say that $\log_{\mathcal{M}}(\pi \mathfrak{m}_K)=2 \mathfrak{m}_K$.
		The following lemma involves more efficient method to show $\log_{\mathcal{M}}(\pi \mathfrak{m}_K)=2 \mathfrak{m}_K$: 
		\begin{lem} \label{l2.9}
			For the uniformizer $\pi$ in the quadratic extension $K$ of $ \mathbb{Q}_2$, we have $\log_{\mathcal{M}}(\pi \mathfrak{m}_K)=2 \mathfrak{m}_K$.
		\end{lem}
		\begin{proof}
			We do know that $\log_{\mathcal{M}}(2 \mathfrak{m}_K)=2 \mathfrak{m}_K$, by Lemma \ref{l2.8}. An element of $\pi \mathfrak{m}_K$ either has $v$-value $1$ or greater; in the latter case, it will be in $2 \mathfrak{m}_K$. So, let's take $z \in \pi \mathfrak{m}_K$, in other words, $v(z) \geq 1$. We need only consider the case of equality, $v(z)=1$. Now, as above, $v(z/2)=0$, so that $z/2 \equiv 1~(\text{mod} \ \mathfrak{m}_K)$. (Here we use the fact that the residue-class field is $\mathbb{F}_2=\mathbb{Z}/2 \mathbb{Z}$.) And of course $-2/2=-1 \equiv 1~(\text{mod} \ \mathfrak{m}_K)$. Thus $(z-2)/2 \equiv 0~ (\text{mod} \  \mathfrak{m}_K)$, and $z-2 \equiv 0~(\text{mod} \ 2 \mathfrak{m}_K)$. The last congruence says $v(z-2) \geq \frac{3}{2}$. Now we get $z+_{\mathcal{M}}(-2)=z-2-2z$, in which $V_2(z-2) \geq \frac{3}{2}$ and $v(-2z) \geq 2$, so that $v(z+_{\mathcal{M}}(-2)) \geq \frac{3}{2}$. So we observe $$\log_{\mathcal{M}}(z)=\log_{\mathcal{M}}(z+_{\mathcal{M}}(-2))-\log_{\mathcal{M}}(-2)=\log_{\mathcal{M}}(z+_{\mathcal{M}}(-2)) \geq \frac{3}{2},$$ which says that $\log_{\mathcal{M}}(z) \in 2 \mathfrak{m}_K$. So $\log_{\mathcal{M}}(\pi \mathfrak{m}_K)=2 \mathfrak{m}_K$.
		\end{proof}
		In these 6 ramified cases, thus, we need only to determine $\log_{\mathcal{M}} (\pi)$ in order to identify $ \log_{\mathcal{M}} (\mathfrak{m}_K)$ because we already know $\log_{\mathcal{M}}(\pi \mathfrak{m}_K)=2 \mathfrak{m}_K$ by Lemma \ref{l2.9} and since the group index $[\mathfrak{m}_K: \pi \mathfrak{m}_K]$ is equal to $2$ with respect to the group structure furnished by $+_{\mathcal{M}}$, we can look at this as saying that the cosets of $\mathfrak{m}_K$ modulo $\pi \mathfrak{m}_K$ are two in number, namely $\pi \mathfrak{m}_K$ and $\pi+_{\mathcal{M}} \pi \mathfrak{m}_K$ and when we apply $\log_{\mathcal{M}}$, these two cosets go over to $\log_{\mathcal{M}}(\pi \mathfrak{m}_K)=2 \mathfrak{m}_K$ and $\log_{\mathcal{M}}(\pi+_{\mathcal{M}} \pi \mathfrak{m}_K)=\log_{\mathcal{M}}(\pi)+\log_{\mathcal{M}}(\pi \mathfrak{m}_K)=\log_{\mathcal{M}}(\pi)+2 \mathfrak{m}_K$.
		\begin{lem} \label{l2.10}
			Let $d$ be one of $\{\pm 6, \pm 2 \}$, and $\pi=\sqrt{d}$. Then $\log_{\mathcal{M}} (\pi)$ is one of the form $u \pi+2 z$, with a unit of $\mathbb{Z}_2$ and $ z \in \mathbb{Z}_2$. Further, in the case $d=2$, $z=0$.
		\end{lem}
		\begin{proof}
			In most cases, we can tell $v_p(\log_{\mathcal{M}}(w))$ from the knowledge of $v_p(w)$, unless $v_p(w)$ is the $\xi$-coordinate of a vertex of the Newton copolygon of $\log_{\mathcal{M}}(x)$. The fact that the $\xi$-coordinates of the vertices of the copolygon, in the case $p=2$, are all at the numbers $\frac{1}{2^n}$ for $n \geq 0$.
			In each case under consideration here, we see that the terms $\frac{\pi^2}{2}$ and $\frac{\pi^4}{4}$ have $ v$ value equal to $0$. This means that we can not ascertain the value of $\log_{\mathcal{M}} (\pi)$ without some direct calculation. We see, however, that the first term in the series $\log_{\mathcal{M}}(\pi)$, namely, $\pi$, dominates. This demonstrates the first claim. \\ For the special situation $d=2$, we use Galois theory, looking at the $\mathbb{Q}_2$-automorphism of $k=\mathbb{Q}_2(\sqrt{2})$, which sends $\pi=\sqrt{2}$ to $-\pi$. We denote this automorphism by $w \to \bar w$. Because the coefficients of $\log_{\mathcal{M}}$ all are in $\mathbb{Q}_2$, we have the relation $\log_{\mathcal{M}}(\bar w)=\overline{\log_{\mathcal{M}}(w)}$ whenever $w \in \mathfrak{m}_K$. In particular, $\log_{\mathcal{M}} (-\pi)=\overline{\log_{\mathcal{M}}(\pi)}$. But we claim that $\pi+_{\mathcal{M}} (-\pi)=-2$: indeed $$ \pi+_{\mathcal{M}} (-\pi)=\mathcal{M}(\pi,-\pi)=\pi+(-\pi)+\pi(-\pi)=-\pi^2=-2.$$   
			By the homomorphic property of $\log_{\mathcal{M}}$, we get $\log_{\mathcal{M}}(\pi)+\log_{\mathcal{M}}(-\pi)=\log_{\mathcal{M}}(-2)=0$, leading to the conclusion that $\log_{\mathcal{M}} (-\pi)=-\log_{\mathcal{M}}(\pi)$. For an element $\alpha+\beta \sqrt{2}$ with $\alpha, \beta \in \mathbb{Q}_2$, when $\bar w=-w$ so that $\alpha-\beta \sqrt 2=-(\alpha+\beta \sqrt{2})$, we conclude $\alpha=0$, since $\{1, \sqrt{2}\}$ is a $\mathbb{Q}_2$-basis of $K$. Or,
			using $\log_{\mathcal{M}}(\pi)=u\pi+2z$ with both $u$ and $z$ in $\mathbb{Z}_2$, since the conjugate of $\log_{\mathcal{M}}(\pi)$ is $-\log_{\mathcal{M}}(\pi)$
			in this case, the $2$-component has to be zero.
		\end{proof}
		\begin{lem} \label{l 4.11}
			Let $d=3$, $\pi=-1+\sqrt{3}$. Then $v_2(\log_{\mathcal{M}}(\pi))=1$, and $\log_{\mathcal{M}}(\pi) \in \mathbb{Z}_2$.
		\end{lem}
		\begin{proof}
			The evaluation of $v_2(\log_{\mathcal{M}} (\pi))$ is computational, and depends only on the polynomial $\pi-\frac{\pi^2}{2}-\frac{\pi^4}{4}-\frac{\pi^8}{8}$ in the evaluation of the $\log_{\mathcal{M}}$, since all others have $v_2$-value $\geq \frac{3}{2}$. \\
			Using the relation $\pi^2+2 \pi-2=0$, we have $\pi^4=4-8 \pi+4 \pi^2=12-16 \pi$. Thus $\pi-\frac{\pi^2}{2}-\pi^4-4=\pi-1+\pi-3+4 \pi=-4+6 \pi$, which has $v_2$-value $\frac{3}{2}$. Since $v_2(\frac{\pi^8}{8})=1$, this term dominates, and the $v$-value of $\log_{\mathcal{M}}(\pi)$ is $1$. \\ 
			The second claim follows by looking at the elements of $1+\mathfrak{m}_K$ corresponding to $\pi$ and $\bar \pi$. These are $1+\pi=\sqrt{3}$ and $-\sqrt{3}$, whose quotient is $-1$, a number whose $p$-adic logarithm is $0$. Thus $\pi-_{\mathcal{M}} \bar \pi=-2$, and when we take logarithms, we get $\log_{\mathcal{M}}(\pi)+\log_{\mathcal{M}}(\bar \pi)=\log_{\mathcal{M}}(\pi)+\overline {\log_{\mathcal{M}}(\pi)}=0,$ so that $\log_{\mathcal{M}} (\pi)$ is its own conjugate and hence in $\mathbb{Q}_2$, the fixed field of conjugation. Since $v_2(\log_{\mathcal{M}}(\pi))=1$, and $\log_{\mathcal{M}}(\pi) \in \mathbb{Q}_2$, we conclude $\log_{\mathcal{M}}(\pi) \in \mathbb{Z}_2$.
		\end{proof}
		\begin{lem} \label{l2.12}
			Let $d=-1$, $\pi=-1+i$. Then $\log_{\mathcal{M}}(\pi)=0$. For $d=-1$, $\log_{\mathcal{M}}(\mathfrak{m}_K)$ has basis equal to the basis of $2 \mathfrak{m}_K$.
		\end{lem}	
		\begin{proof}
			The element of $1+\mathfrak{m}_K$ corresponding to $\pi$ is $1+\pi=i$, multiplicatively of order $4$, so that its $\log_{\mathcal{M}}(\pi)=\log_{\mathcal{M}}(-1+i)=\log_{\mathcal{M}}(-1)+\log_{\mathcal{M}}(i)=0$. In another way, since $\pi=-1+i$, we have $\pi^2+2 \pi=2$. Now, $2\log_{\mathcal{M}}(\pi)=\log_{\mathcal{M}}(\pi+_{\mathcal{M}}\pi)=\log_{\mathcal{M}}(\pi+\pi+\pi^2)=\log_{\mathcal{M}}(2 \pi+\pi^2)=\log_{\mathcal{M}}(-2)=0$, by Lemma \ref{l2.1}.
			Next, the conclusion is just another way of saying that $\log_{\mathcal{M}}(\mathfrak{m}_K)=2 \mathfrak{m}_K$. we know that $\log_{\mathcal{M}}(\pi \mathfrak{m}_K)=2 \mathfrak{m}_K$. Consider $z \in \mathfrak{m}_K$ now, then as before, $v(z/\pi) \equiv 1~(\text{mod}~ \mathfrak{m}_K)$, and thus $z/\pi+1 \equiv 0~(\text{mod} ~\mathfrak{m}_K)$, again using the fact that the residue field is $\mathbb{F}_2$. As a result, $z+\pi \equiv 0~(\text{mod}~ \pi \mathfrak{m}_K)$. Again we calculate $z+_{\mathcal{M}} \pi=z+\pi+\pi z$ in which both $z+\pi$ and $\pi z$ are in $\pi \mathfrak{m}_K$, so that $z+_{\mathcal{M}} \pi \in \pi \mathfrak{m}_K$.Now we see that: $$\log_{\mathcal{M}}(z)=\log_{\mathcal{M}}(z+_{\mathcal{M}} \pi)-\log_{\mathcal{M}}(\pi)=\log_{\mathcal{M}}(z+_{\mathcal{M}} \pi) \in \log_{\mathcal{M}}(\pi \mathfrak{m}_K)=2 \mathfrak{m}_K.$$ Hence $\log_{\mathcal{M}}(\mathfrak{m}_K)$ has basis equal to the basis of $2 \mathfrak{m}_K$ for $d=-1$.  
		\end{proof}
		
		We are now ready to describe $\log_{\mathcal{M}}(\mathfrak{m}_K)$ in the $6$ ramified cases $d=-1,\pm 2, 3, \pm 6$. We would describe the set $\log_{\mathcal{M}}(\mathfrak{m}_K)$ by giving a basis $\{\alpha, \beta\}$, because it is a free $\mathbb{Z}_2$-module of rank $2$. We know that $\log_{\mathcal{M}}(\mathfrak{m}_K)$ is a group containing $\log_{\mathcal{M}}(\pi \mathfrak{m}_K)=\log_{\mathcal{M}}(2 \mathfrak{m}_K)=2\mathfrak{m}_K,$ with the associated index being $2$. Here we note that in ramified extensions, we have a $\mathbb{Z}_2$-basis $\{1, \pi \}$ for $\mathcal{O}_K$ and to get a $\mathbb{Z}_2$-basis of $\pi^m \mathcal{O}_K$, we just take $\{\pi^m, \pi^{m+1}\}$. In case $m=2r$ is even, an equally good basis will be $\{2^r,2^r \pi \}$, and in case $m=2r+1$ is odd, an equally good basis will be $\{2^r \pi, 2^{r+1}\}$. We will use the basis $\{2 \pi,4\}$ for $2\mathfrak{m}_K$ to see how it extends when we adjoin $\log_{\mathcal{M}}(\pi)$.
		\begin{thm} \label{t2.13}
			If $K=\mathbb{Q}_2(\sqrt{d})$ for $d=-1, \pm 2,3, \pm 6$, where $\pi=\sqrt d$ in case of even $d$ and $\pi=-1+\sqrt d$ for odd $d$, then: 
			\begin{enumerate}
				\item For $d=-1$, $\log_{\mathcal{M}}(\mathfrak{m}_K)$ has basis $\{2 \pi, 4\}$.
				\item For $d=-2$, $\log_{\mathcal{M}}(\mathfrak{m}_K)$ has basis $\{ \pi u+2z,4 \}$ with $z \in \mathbb{Z}_2$.
				\item For $d=2$, $\log_{\mathcal{M}}(\mathfrak{m}_K)$ has basis $\{\pi,
				4\}$.
				\item For $d=3$, $\log_{\mathcal{M}}(\mathfrak{m}_K)$ has basis $\{2 \pi,2\}$.
				\item For $d=-6$, $\log_{\mathcal{M}}(\mathfrak{m}_K)$ has basis $\{\pi+2z,4\}$ with $z \in \mathbb{Z}_2$.
				\item For $d=6$, $\log_{\mathcal{M}}(\mathfrak{m}_K)$ has basis $\{\pi+2z, 4\}$, with $z \in \mathbb{Z}_2$.
			\end{enumerate}
			
		\end{thm}
		\begin{proof}
			\begin{enumerate}
				
				\item[(1)] The basis of $2 \mathfrak{m}_K$ is $\{2 \pi, 4\}$. By Lemma \ref{l2.12}, the basis of $\log_{\mathcal{M}}(\mathfrak{m}_K)$ is equal to the basis of $2 \mathfrak{m}_K$. Therefore, the basis of $\log_{\mathcal{M}}(\mathfrak{m}_K)$ is $\{2 \pi, 4\}$.
				
				Alternatively, we can directly show that $\log_{\mathcal{M}}(\mathfrak{m}_K)=2 \mathfrak{m}_K$. For, we know that from Lemma \ref{l2.9}, $\log_{\mathcal{M}}(\pi \mathfrak{m}_K)=2 \mathfrak{m}_K$. Consider $z \in \mathfrak{m}_K$, then $v(z/\pi) \geq 0$, and as before the only one case that concerns us is when $v(z/\pi)=0$ so that $z/\pi \equiv 1$ (mod $\mathfrak{m}_K$), and thus $z/\pi+1 \equiv 0$ (mod $\mathfrak{m}_K$), by using the fact that $1=-1$ (mod 2) because  $\mathbb{Q}_2(\sqrt{-1})$ is ramified extension of degree 2 and its residue field is $\mathbb{F}_2$. As a result, $z+\pi \equiv 0$ (mod $\pi \mathfrak{m}_K$). Again we calculate $z+_{\mathcal{M}} \pi=\mathcal{M}(z,\pi)=z+\pi+\pi z$, in which both $z+\pi$ and $\pi z$ are in $\pi \mathfrak{m}_K$, so that $z+_{\mathcal{M}} \pi \in \pi \mathfrak{m}_K$. Finally, \begin{align*}
					\log_{\mathcal{M}}(z)&=\log_{\mathcal{M}}(z+_{\mathcal{M}}\pi)-\log_{\mathcal{M}}(\pi) \\
					&=\log_{\mathcal{M}}(z+_{\mathcal{M}}\pi) \in \log_{\mathcal{M}}(\pi \mathfrak{m}_K)=2 \mathfrak{m}_K,~(\log_{\mathcal{M}}(\pi)=0)
				\end{align*}
				\item[(2)] This follows from Lemma \ref{l2.10}.
				\item[(3)] We have already seen, in case $K=\mathbb{Q}_2(\sqrt{2})$, that $\log_{\mathcal{M}}(\pi)=u \pi$ where $u$ is a $\mathbb{Z}_2$-unit. We also know that $\log_{\mathcal{M}}(\pi \mathfrak{m}_K)=2 \mathfrak{m}_K$, which has basis $\{2 \pi,4\}$. Since the \enquote{new} element of $\log_{\mathcal{M}}(\mathfrak{m}_K)$ is $\log_{\mathcal{M}}(\pi)=u \pi$, there is also an element of $\mathfrak{m}_K$, namely $w=[u^{-1}]_{\mathcal{M}} (u \pi)$ whose $\log_{\mathcal{M}}(w)=\pi=\sqrt{2}$. Recall that $[u^{-1}]_{\mathcal{M}}$ is the multiplication-by-$u^{-1}$ map, which is an endomorphism of the 1-dimensional multiplicative formal group $\mathcal{M}$. Also $\log_{\mathcal{M}}([2]_{\mathcal{M}}(w))=2 \pi$, which belong to $\log_{\mathcal{M}}(\mathfrak{m}_K)$. This demonstrate the basis $\{\pi,4 \}$.
				\item[(4)] For the case $d=3$, $\pi=\sqrt 3-1$. We know $\log_{\mathcal{M}}(\pi \mathfrak{m}_K)=2 \mathfrak{m}_K$ by Lemma \ref{l2.9}, which has basis $\{2 \pi,4 \}$.  The image of $\log_{\mathcal{M}}$ contains at least $2\pi \mathcal{O}_K=2\mathfrak{m}_K$, and the important thing to note that on this set, the 2-adic logarithm and 2-adic exponential are perfect inverses of each other. Now, the index of the quotient group $\mathfrak{m}_K/2\mathfrak{m}_K$ is $4$ by Lemma \ref{l2.7}. Note that the quotient group $\mathfrak{m}_K/2\mathfrak{m}_K$ is $\mathbb{F}_2$-vector space with basis $\{\pi, 2\}$. Let us now look at what the logarithm $\log_{\mathcal{M}}$ does to the elements of $\mathfrak{m}_K/2\mathfrak{m}_K$. We have already seen that $\log_{\mathcal{M}}(\pi)=2u$ for a unit $u\in\mathbb{Z}_2$, and thus $\log_{\mathcal{M}}([1/u](\pi)=2$. So $2$ is in the image of $\log_{\mathcal{M}}$. By direct calculation, one sees that $\log_{\mathcal{M}}(-\pi)$ has the form $2z + 2w\pi$, where $z$ is something in $\mathbb{Z}_2$ and $w$ is a unit in $\mathbb{Z}_2$. And that gives us the \enquote{new} element of the image of the logarithm $\log_{\mathcal{M}}$, giving as basis of $\log_{\mathcal{M}}(\mathfrak{m}_K)$ the pair $\{2 \pi, 2\}$.
				\item[(5)] The proof follows from Lemma \ref{l2.10}.
				\item[(6)] The proof follows from Lemma \ref{l2.10}.
			\end{enumerate} 
		\end{proof} 
		It only remains to dispose of the single unramified case, $K=\mathbb{Q}_2(\sqrt{-3})$ which is proved in the following theorem:
		\begin{thm} \label{t2.14}
			A $\mathbb{Z}_2$-basis of $\log_{\mathcal{M}}(\mathfrak{m}_K)$ is $\{2, 4 \omega\}$ in the unramified case $K=\mathbb{Q}_2(\omega)$, where $\omega^2+\omega+1=0$. 
		\end{thm}
		\begin{proof}
			Here, $2$ is a prime, but the residue field is $\mathbb{F}_4=\mathbb{F}_2(\omega)$ with $\omega^2+\omega+1=0$. Back here in characteristic zero, there is still an $\omega \in \mathcal{O}_K$, $\omega^2+\omega+1=0$, and we may take $\omega=\frac{-1+\sqrt{-3}}{2}$, a cube root of unity. \\
			We still have the result that $\log_{\mathcal{M}}(2\mathfrak{m}_K)=2\mathfrak{m}_K$, but since $\mathcal{O}_K$ has basis $\{1, \omega \}$ and $\mathfrak{m}_K=2 \omega$, the module $2\mathfrak{m}_K$ has $\mathbb{Z}_2$-basis $\{4,4 \omega\}$. Since $\mathfrak{m}_K$ has basis $\{2,2 \omega\}$, we have $[\mathfrak{m}_K:2 \mathfrak{m}_K]=4$. Let us use our knowledge that $\log_{\mathcal{M}}(-2)=0$ to say that $[\log_{\mathcal{M}}(\mathfrak{m}_K):\log_{\mathcal{M}}(2\mathfrak{m}_K)] \leq 2$. The only question is what $\log_{\mathcal{M}}(2 \omega)$ might be. We get \begin{align*}
				\log_{\mathcal{M}}(2 \omega) &=2 \omega-\frac{4 \omega^2}{2}+\frac{8 \omega^3}{3}-\frac{16 \omega^4}{4}+\cdots \\ &=2 \omega-2 \omega^2+\frac{8 \omega^3}{3}-\frac{16 \omega^4}{4}+\cdots,
			\end{align*}
			where the missing terms all have $v_2$-value greater than $1$. Indeed, the only terms that interest us are $2 \omega-2 \omega^2$. But $\omega^2=-\omega-1$, the $v$-value of $2+4 \omega$ is $1$, so that $\log_{\mathcal{M}}(\mathfrak{m}_K)$ in this case has basis $\{2,4 \omega\}$. 
		\end{proof}
		
		\section{Applications} \label{s3}
		The basis of vector space is as important as the basis of the free module. So Theorems \ref{t2.13} and \ref{t2.14} give important directions. We write a few applications of our results as follows:

		Let $p$ be a prime number, $\mathbb{Q}_p$ the field of $p$-adic numbers, and $\mathbb{C}_p$ the completion
		of the algebraic closure of $\mathbb{Q}_p$. Let $U_p$ be the units $(1+\mathfrak{m}_{\mathbb{C}_p})$ of $\mathbb{C}_p$, the $p$-adic version of Baker theorem is:
		\begin{thm} \label{t1.1} (Brumer)
			Let $\alpha_1, \cdots, \alpha_n$ be elements of principal units $U_p$ which are algebraic over the rationals $\mathbb{Q}$ and their $p$-adic logarithms $\log_p(\alpha_1), \cdots, \log_p(\alpha_n)$ are linearly independent over $\mathbb{Q}$. These logarithms $\log_p(\alpha_1), \cdots, \log_p(\alpha_n)$ are
			then linearly independent over the algebraic closure $A$ of $\mathbb{Q}$ in $\mathbb{C}_p$.
		\end{thm}
		
		We prove a thin result similar to the above result for $p=2,~n=2$.
		\begin{thm} \label{t3.2}
			If $\alpha_1, \alpha_2$ are two elements of the group of principal units $U_1=1+\mathfrak{m}_K$ of $K=\mathbb{Q}_2(\sqrt{-1})$ such that their 2-adic logarithms $\log_2(\alpha_1), \log_2(\alpha_2)$ are linearly independent over $\mathbb{Q}$. Then their logarithms $\log_2(\alpha_1), \log_2(\alpha_2)$ are
			linearly independent over the algebraic closure $A$ of $\mathbb{Q}$ in $\mathbb{Q}_2$.
		\end{thm}
		\begin{proof}
			Note the uniformizer of $K=\mathbb{Q}_2(\sqrt{-1})$ was $\pi=-1+\sqrt{-1}$. Since $\{2 \pi,4 \}$ is a $\mathbb{Z}_2$-basis of $\log_2(U_1)$ by Theorem \ref{t2.13}, we can write $$\log_2(\alpha_1)=2 \pi a_1+4b_1,~~\log_2(\alpha_2)=2\pi a_2 +4b_2,~~a_i,b_i \in \mathbb Z_2.$$
			But given that $\log_2(\alpha_1), \log_2(\alpha_2)$ are linearly independent over $\mathbb{Q}$, therefore
			\begin{align} \label{e1}
				\det \begin{pmatrix} a_1 & b_1 \\ a_2 & b_2 \end{pmatrix} \neq 0
			\end{align}
			
			We claim $\log_2(\alpha_1), \log_2(\alpha_2)$ are also linearly independent only over $\mathbb Q_2 \cap \bar{\mathbb{Q}}$. For,
			if $\beta_1, \beta_2 \in \mathbb Q_2 \cap \bar{\mathbb{Q}}$ such that 
			\begin{align} \nonumber&\beta_1 \log_2(\alpha_1)+\beta_2\log_2(\alpha_2)=0 \\ \Rightarrow \nonumber&\beta_1(2 \pi a_1+4b_1)+\beta_2(2 \pi a_2+b_2)=0 \\ \Rightarrow \nonumber&4(b_1\beta_1+b_2\beta_2)+2 \pi(a_1\beta_1+a_2\beta_2)=0,~a_i,b_i \in \mathbb Z_2~\text{and}~\beta_1, \beta_2 \in \mathbb Q_2 \cap \bar{\mathbb{Q}} \\ \Rightarrow \nonumber & b+\pi a=0,~~\text{where}~ b=4(b_1\beta_1+b_2\beta_2),~a=2 (a_1\beta_1+a_2\beta_2) \in \mathbb Q_2 \\ \Rightarrow \nonumber&b=a=0,~\text{since}~\pi \notin \mathbb Q_2. \\ \Rightarrow  & \nonumber \left\{\begin{aligned} a_1\beta_1+a_2\beta_2=0 \\ b_1\beta_1+b_2\beta_2=0  \end{aligned}\right.
			\end{align}
			By relation \eqref{e1}, this system has a trivial solution $\beta_1=\beta_2=0$.
			Therefore $\log_2(\alpha_1), \log_2(\alpha_2)$ are linearly independent over the algebraic closure $A$ of $\mathbb{Q}$ in $\mathbb{Q}_2$. This proves the theorem.
		\end{proof}
		\begin{rem}
			In the above theorem, we can replace $\mathbb{Q}_2(\sqrt{-1})$ by other 6 quadratic extensions of $\mathbb{Q}_2$.
		\end{rem}

		\begin{prop}
			The $2$-adic logarithm induces a measure on $4 \mathbb{Z}_2$ from a measure on $1+2 \mathbb{Z}_2$.
			
		\end{prop}
		\begin{proof} By \cite[Sec.~12.2,~Chap.~12]{WA}, if $h: X \to Y$ is continuous and $d \phi$ is a measure on $X$, then we obtain a measure $d \psi$ on $Y$ by defining 
			\begin{align} \label{e9}
				\int_Yf d \psi=\int_X f(h(X)) d \phi,
			\end{align}
			for some measurable function $f$. 
			
			Take $X=1+2 \mathbb{Z}_2,~Y=4\mathbb{Z}_2$, and $h=\log_2$. Recall that $p$-adic logarithm is a continuous function. Infact, the 2-adic logarithm $\log_2: 1+2 \mathbb{Z}_2 \to 4 \mathbb{Z}_2$ is continuous in Lemma \ref{l2.8}. Thus from \eqref{e9}, we have 
			\begin{align*}
				\int_{4 \mathbb{Z}_2} f d \psi&=\int_{1+2 \mathbb{Z}_2}f(\log_2(1+2 \mathbb{Z}_2))d \phi, \\
				\Rightarrow \int_{4 \mathbb{Z}_2} f d \psi&=\int_{1+2 \mathbb{Z}_2}f(4 \mathbb{Z}_2)d \phi,~(\text{by Lemma \ref{l2.8}}).
			\end{align*}
			Therefore, we obtain a measure on $4 \mathbb{Z}_2$ induced from a 2-adic measure on $1+2 \mathbb{Z}_2$.
		\end{proof}

		\begin{prop}
			 For two arbitrary primes $p_1,p_2 \equiv 5$ (mod 8), there is a 2-adic integer $z$ such that $p_1=p_2^z$.
		\end{prop}
		\begin{proof}
			For, $p_1 \equiv 5$ (mod 8), we have 
			\begin{align*}
				v_2(\log_2(p_1))&=v_2(\log_2(p_1)-\log_2(1)),~\text{since}~\log_2(1)=0 \\
				&=v_2(p_1-1),~\text{since $\log_2$ is isometry from $1+4 \mathbb{Z}_2$ to $4 \mathbb{Z}_2$ by Lemma \ref{l2.8}} \\
				&\equiv v_2(5-1)~\text{mod 8} \equiv 2~\text{mod 8}.
			\end{align*}
			Similary, $v_2(\log_2(p_2)) \equiv 2~\text{mod 8}$. Thus $v_2\left(\frac{\log_2(p_1)}{\log_2(p_2)}\right)=v_2(1)=0$, and so $\frac{\log_2(p_1)}{\log_2(p_2)} \in \mathbb{Z}_2$.

			We know $\log_2(x^n)=n \log_2(x)$ holds for $n \in \mathbb{Z}$, $x \in 1+4 \mathbb{Z}_2$. Now the map $(x,n) \mapsto x^n$ from $(1+4 \mathbb{Z}_2) \times \mathbb{Z} \to 1+4 \mathbb{Z}_2$ is continuous with respect to $2$-adic topology, and thus extends to $(1+4 \mathbb{Z}_2) \times \mathbb{Z}_2 \to 1+4 \mathbb{Z}_2$. Composing with continuous map $\log_2$, we get $\log_2(x^z)=z \log_2(x)$ for all $z \in \mathbb{Z}_2$. So taking the 2-adic logarithm both sides of $p_1=p_2^z$, we obtain $$\log_2(p_1)=z \log_2(p_2) \Rightarrow z=\frac{\log_2(p_1)}{\log_2(p_2)} \in \mathbb{Z}_2.$$ 
			Thus the equation $p_1=p_2^z$ has a solution in $\mathbb{Z}_2$ i.e., $z$ is a 2-adic integer satisfying $p_1=p_2^z$.
		\end{proof}
		
		\subsection*{$p$-adic logarithm in Iwasawa theory}
		
		We repeat from the introduction section that, the reference to Iwasawa is the most important application of $p$-adic logarithm. In \cite{KI}, Iwasawa has extensively studied certain types of modules made up of the images of $p$-adic logarithm on the principal units $U_n$ of the field $\mathbb{Q}_p(\zeta_{p^{n+1}})$ and used those special modules in \cite{KI1} to prove explicit formulas for the norm residue symbol, e.g., \cite[Lemma~1]{KI}, \cite[Theorem~1]{KI1}. There are other similar results along the lines of Iwasawa in \cite[Theorem~1.10]{BT}, \cite[Section~13.8]{WA} and in \cite{JC}, where the image of $p$-adic logarithm on the principal units $U_n$ of the cyclotomic extension of $\mathbb{Q}_p$ has played a crucial role.

		\subsection*{Conclusion} The paper takes an initial effort to compute the image of $p$-adic logarithm on the principal units of quadratic extensions of $\mathbb{Q}_2$. We find the $\mathbb{Z}_2$-bases of $\log_p(1+\mathfrak{m}_K)$. In the later work \cite{MA}, the authors explicitly compute the images of $\log_p(1+\mathfrak{m}_K)$ on the principal units in all 7 quadratic extensions of $\mathbb{Q}_2$ as well as for cyclotomic extensions of $\mathbb{Q}_p$ for all $p$. Iwasawa's work is our main motivaton to study the image of $p$-adic logarithm on the group of principal units.
		\subsection*{Question} It would be interesting to compute the bases of $p$-adic logarithm on the principal units of any finite extensions of $\mathbb{Q}_p$. 
		\subsection*{Acknowledgement} The authors are deeply grateful to Professor Jonathan Lubin for helpful correspondence and suggestions. The first author is grateful to \textit{CSIR}, Govt. of India, for the grant with File no.-09/025(0249)/2018-EMR-I.
	

\begin{thebibliography}{50}
		\bibitem{MA} M. A. Sarkar and A. A. Shaikh, On the image of $p$-adic logarithm, https://arxiv.org/pdf/1904.09850.pdf.
		\bibitem{NK} N. Koblitz, \text{p-adic Numbers, p-adic Analysis, and Zeta-functions}, Springer, New York, 1984.
		\bibitem{AR} A. M. Robert, A course in $p$-adic Analysis, Springer, New York, 2000.
		\bibitem{LB} L. Berger, An introduction to the theory of $p$-adic representations, Lecture series in \enquote{Dwork trimester}(2001), Padova, https://arxiv.org/abs/math/0210184. 
		\bibitem{M}  M. Hazewinkel, \textit{Formal Groups and Applications}, AMS Chelsea Publishing, 1978.
		\bibitem{J} J. H. Silverman, \textit{The Arithmetic of Elliptic curves}, Springer Publication, 1992.
		\bibitem{GS} G. Shimura, Arithmetic of Quadratic Forms, Springer, New York, 2009.
		\bibitem{BT} B. Angles and T. Herreng, On a result of Iwasawa, International Journal of Number Theory, 2007.
		\bibitem{FY} C. Fieker and Y. Zhang, An application of the-adic analytic class number formula, LMS Journal of Computation and Mathematics 19(1), 217-228, 2016.	
		\bibitem{AB} A. Brumer, On the units of algebraic number fields, Mathematika 14(2), 121-124, 1967. 
		\bibitem{KI1} K. Iwasawa, On some modules in the theory of cyclotomic fields, J. Math. Soc. Japan 16, 42-82, 1964.
		\bibitem{KI} K. Iwasawa, On explicit formulas for the norm residue symbol, J. Math. Soc. Japan 20, 151-165, 1968.
		
		\bibitem{JC} J. Coates, $p$-adic $L$-functions and Iwasawa's theory, in Algebraic number fields (Durham symposium, 1975, ed. by A. Fr$\ddot{o}$hlich), 269-353, Academic Press, London 1977.
		\bibitem{WA} L.C. Washington, Introduction to cyclotomic fields, Sprinjer Verlag, 1997.
		\bibitem{GG} G. Georges, The p-adic Kummer-Leopoldt constant-Normalized
		p-adic regulator, International Journal of Number Theory, 14(2), 329-337, 2018.
		\bibitem{BM} B. Mazur, W. Stein and J. Tate, Computation of $p$-adic Heights and Log Convergence, https://people.math.harvard.edu/~mazur/preprints/pheight.pdf. 

		
		
		
	\end{thebibliography}
\end{document}